\documentclass{amsart}

\usepackage{amsmath, amssymb, amsfonts, amscd, amsthm, mathrsfs}
\usepackage{tikz}
\usetikzlibrary{arrows, automata, matrix, positioning}
\usepackage{graphicx}
\usepackage{float}
\usepackage{epstopdf}
\usepackage[all]{xy}

\usepackage{url}
\makeatletter
\def\url@leostyle{%
  \@ifundefined{selectfont}{\def\UrlFont{\sf}}{\def\UrlFont{\small\ttfamily}}}
\makeatother
\newtheorem{theorem}{Theorem}[section]
\newtheorem{lemma}[theorem]{Lemma}
\newtheorem{proposition}[theorem]{Proposition}
\newtheorem{corollary}[theorem]{Corollary}
\newtheorem{fact}[theorem]{Fact}

\newtheorem{question}[theorem]{Question}

\newtheorem{definition}[theorem]{Definition}

\newtheorem{remark}[theorem]{Remark}
\newtheorem{example}[theorem]{Example}

\newenvironment{cproof}{%
  \par\noindent\textbf{Proof. }\ignorespaces
}{%
  \hfill\qed\\
}

\newtheorem*{theorem*}{Theorem}
\newtheorem*{lemma*}{Lemma}
\newtheorem*{proposition*}{Proposition}
\newtheorem*{corollary*}{Corollary}
\newtheorem*{fact*}{Fact}
\newtheorem*{facts*}{Facts}
\newtheorem*{claim*}{Claim}

\newtheorem*{definition*}{Definition}
\newtheorem*{rmkdf*}{Remark and Definition}
\newtheorem*{remark*}{Remark}
\newtheorem*{example*}{Example}
\newtheorem*{examples*}{Examples}
\newtheorem*{exercise*}{Exercise}

\usepackage{graphicx} 

\newcommand{\R}{\mathbb{R}}

\newcommand{\cM}{\mathcal{M}}

\newcommand{\Lg}{\mathscr{L}}
\newcommand{\Rg}{\mathscr{R}}
\newcommand{\Dg}{\mathscr{D}}
\newcommand{\Jg}{\mathscr{J}}
\newcommand{\Hg}{\mathscr{H}}
\newcommand{\bi}[1]{S^1 #1 S^1}
\newcommand{\es}[1]{S^1 #1}
\newcommand{\di}[1]{#1 S^1}
\newcommand{\M}[4]{\mathscr{M}(#1,#2,#3,#4)}

\usepackage[colorlinks=true, linkcolor=blue, citecolor=blue, urlcolor=blue]{hyperref}

\title{On Definably Compact Semigroups in o-Minimal Structures}
\author{Eduardo Magalhães}
\date{May 2025}

\address{Department of Mathematics, University of Porto, Portugal}
\email{eduardomag79@gmail.com}

\begin{document}

\begin{abstract}
In this paper, we study definably compact semigroups in o-minimal structures, aiming to extend the theory of definable groups to a broader algebraic setting. We show that any definably compact semigroup contains idempotents and admits a unique minimal ideal as well as minimal left and right ideals, all of which are definable and definably compact. We prove that cancellativity characterizes definable groups among definably compact semigroups, and that any definably compact subsemigroup of a definable group is itself a subgroup. In the completely simple case, we classify definable subsemigroups via a definable Rees matrix decomposition and analyse when results from the theory of definable groups extend to this setting.
\end{abstract}

\maketitle

\section{Introduction}

In 1986, the definition of an o-minimal structure was introduced  by A. Pillay and C. Steinhorn in ``Definable Sets in Ordered Structures I'' \cite{pillay_definable_1986}, motivated by the work of L. van den Dries in \cite{den_dries_remarks_1984}. Two years later, in 1988, A. Pillay published ``On Groups and Fields Definable in o-Minimal Structures'' \cite{pillay_groups_1988}, which is widely regarded as the beginning of the study of groups definable in o-minimal structures. 

Since then, the theory of definable groups in o-minimal structures has evolved into a robust and elegant framework that highlights how o-minimality restricts and influences the algebraic behaviour of definable groups. This success prompts a natural question: to what extent can these results and techniques be extended to more general algebraic objects, such as semigroups? Surprisingly, despite the ubiquity of semigroups in algebra and topology, their definable counterparts in o-minimal structures appear to be largely undeveloped: to our knowledge, definable semigroups have not previously been studied in the o-minimal setting in a systematic way. This paper aims to begin filling that gap.

In this paper, we investigate whether classical theorems concerning semigroups and topological semigroups remain valid when translated to the definable setting within o-minimal structures, and whether results about definable groups in o-minimal structure carry over to definable semigroups. In doing so, we not only aim to generalize the theory of definable groups, but also contribute to the broader program of understanding algebraic and topological behaviour under the geometric constraints of o-minimality.

While semigroups lack inverses and many structural tools available for groups, their study has long revealed subtle algebraic and combinatorial behaviour, especially with the addition of extra structure, like a topology compatible with the semigroup operation. In particular, compactness plays a central role in the classical theory of topological semigroups. Motivated by this, the focus of this paper is the study of definably compact definable semigroups.

Throughout this paper, we will be working in a fixed sufficiently saturated o-minimal structure $\cM$, and all definable semigroups and groups are assumed to carry a Hausdorff definable space structure compatible with their operations. As such, a key feature of our approach is that we work with definable space structures on all semigroups and, in particular, on all groups. This is a significant generalization of the usual setting in which definable groups are studied, where one typically works with a definable manifold structure making the group a topological group (as in Pillay’s topology). By contrast, we allow general Hausdorff definable spaces.

We start by reviewing some preliminary definitions and results in Section 2 that we will use heavily in the rest of this paper, namely the definition of a definable space, definable manifold and definable compactness. Section 2 also contains a quick introduction to the semigroup theoretic definitions and results that will be used. 

In section 3 we begin the study of definably compact semigroups. In parallel with classical results concerning compact semigroups, we prove:

\begin{theorem*}[\ref{def compact semigroup has idempotent}]
    Let $S$ be a definably compact definable semigroup. Then $S$ has at least one  idempotent.
\end{theorem*}
\begin{theorem*}[\ref{unique minimal ideal}]
    Let $S$ be a definably compact definable semigroup. Then $S$ has a unique minimal ideal. Furthermore, that minimal ideal is definable and definably compact.
\end{theorem*}

In Section 4, we study definably compact completely simple semigroups, which is the case of the minimal ideal of any definably compact semigroup. In particular, we prove that some properties of definable groups hold in definably compact completely simple semigroups if and only if we impose strong finiteness properties (namely, $\Hg$-finiteness).

\section{Preliminaries}

\subsection{Definable Spaces and Definable Compactness}

We will assume that the reader is familiar with the basics of o-minimality (see \cite{dries_tame_1998} for a good introduction to the subject). Additionally, throughout this paper, definable will always mean definable with parameters, unless stated otherwise.

The notion of a definable space and definable manifold play a central role in this work. As such, although the definitions are well known, we will nonetheless explicitly state them.

\begin{definition*}
     Let $\cM$ be an o-minimal structure, and let $X$ be an arbitrary set (not necessarily definable). A definable atlas on $X$ relative to $\cM$ is a finite family $(U_i, \phi_i)_{i  \in I}$ such that:

    \begin{enumerate}
        \item For each $i \in I$, $U_i\subseteq X$ and $X = \bigcup_{i \in I} U_i$;

        \item For each $i \in I$, $\phi_i(U_i) \subseteq M^{n_i}$ is definable and $\phi_i: U_i \to \phi_i(U_i)$ is a bijection;
        
        \item For each $i,j \in I$, $\phi_i(U_i\cap U_j)$ is definable and open in $\phi_i(U_i)$ with respect to the topology induced by $\cM$;
        
        \item For each $i,j \in I$, each map $\phi_{ji}:= \phi_j \circ \phi_i^{-1}:= \phi_i(U_i\cap U_j) \to  \phi_j(U_i\cap U_j)$ is a definable homeomorphism.

    \end{enumerate}
\end{definition*}

    Given a set $X$ equipped with definable atlas $(U_i, \phi_i)_{i  \in I}$, we can uniquely define a topology on $X$ by: $V \subseteq X$ is open if and only if, for all $i \in I$, the set $\phi_i(U_i \cap V)$ is open in $\phi(U_i)$. We refer to this unique topology as the topology induced by the atlas $(U_i, \phi_i)_{i  \in I}$.

\begin{definition*}
    A definable space is simply a set $X$ equipped with a definable atlas $(U_i, \phi_i)_{i  \in I}$.

    Additionally, if the definable atlas satisfies: for all $i \in I$, $\phi_i(U_i) \subseteq M^n$ is open, for some fixed $n \geq 0$, then we say that $X$ is a definable $n$-manifold, or simply definable manifold. 
\end{definition*}

In particular, note that any definable manifold is also a definable space.

\begin{example}
    Let $\cM = (\R, +, \cdot)$ and let $X = [0,1)$. Let $\{(U_1,\phi_1),(U_2,\phi_2)\}$ be a definable atlas given by:
    \begin{itemize}
        \item $U_1 = \{0\}$ and $\phi_1 : \{0\} \to \{0\}$ is the identity;
        \item $U_2 = (0,1)$ and $\phi_2 : (0,1) \to (0,1)$ is the identity;
    \end{itemize}

    Then $X$ equipped with the definable atlas $(U_i,\phi_i)_{i=1,2}$ is a definable space, but not a definable manifold.

    On the other hand, consider the definable atlas $\{(V_1,\psi_1),(V_2,\psi_2)\}$ given by:
    \begin{itemize}
        \item $V_1 = [0, 0.1) \cup (0.9, 1)$ and $\psi_1: V_1\to (0,0.2)$ given by \begin{align*}
            \psi_1(t) = \begin{cases}
                t + 0.1 & \text{if $t \in [0,0.1)$}\\
                t - 0.9 & \text{if $t \in (0.9,1)$}
            \end{cases}.
        \end{align*}
        \item $V_2 = (0,1)$ and $\psi_2 : (0,1) \to (0,1)$ is the identity;
    \end{itemize}

    Then $X$ equipped with the definable atlas $(V_i,\psi_i)_{i=1,2}$ is a definable $1$-manifold. 
    
    In particular, note that $(V_i,\psi_i)_{i=1,2}$ and $(U_i,\phi_i)_{i=1,2}$ induce different topologies on $X$: $X$ equipped with $(U_i,\phi_i)_{i=1,2}$ has $0$ as an isolated point, while $X$ equipped with $(U_i,\phi_i)_{i=1,2}$ is homeomorphic to $\mathbb{S}^1$ and thus has no isolated points.
\end{example}

Recall that a definable group in $\cM$ is a group $(G,\cdot)$ such that $G\subseteq M^n$ is a definable set, for some $n\geq 1$, and $\cdot :G \times G \to G$ is a definable map. 

In \cite{pillay_groups_1988}, A. Pillay proved that any definable group in an o-minimal structure admits a unique definable manifold topology $t$ (i.e. a topology arising from a definable manifold structure) that makes it a topological group. This remarkable result was essential to many advances in the study of definable groups. Among the results that follow from the existence of such topology, we would like to remark the following.

Recall that we say that $G$ has the descending chain condition (DCC) if any descending chain of definable subgroups eventually stabilizes.

\begin{fact}\label{facts about definable groups}
    Let $G$ be a definable group in some o-minimal structure. Then:
    \begin{enumerate}
        \item Any definable subgroup of $G$ is $t$-closed (i.e. closed with respect to the topology $t$);
        \item A subgroup $H$ is $t$-open if and only if $\dim H = \dim G$;
        \item $G$ has the DCC;
        \item There are only finitely many definable subgroups $H$ of $G$ with $\dim H = \dim G$.
    \end{enumerate}
\end{fact}

 The preceding fact follows from Corollary 2.8 and Lemma 2.10 in \cite{pillay_groups_1988}, Fact 1.13 in \cite{peterzil_definably_2000} and Theorem 2.6 in \cite{strzebonski_euler_1994}, respectively.

\begin{remark}\label{remark about definable}
    The definition of definable space given above is the general one: the underlying set need not itself be a definable subset of some Cartesian power of $M$. In that generality, definability of subsets of a definable space is defined through the atlas (see Chapter 10, \S 1 of \cite{dries_tame_1998}).
    
    From now on, however, every definable space we use will be of the following more concrete form. The underlying set will be a definable set \(X\subseteq M^n\), and the chosen atlas \((U_i,\phi_i)_{i\in I}\) will be definable in the sense that each \(U_i\subseteq X\) and each chart map \(\phi_i\) is definable. For such spaces, the atlas-induced notion of definable subset agrees with the ambient one: a subset \(Y\subseteq X\) is definable if and only if it is definable as a subset of \(M^n\).
\end{remark}

Finally, we will briefly discuss a generalization of compactness to our setting. Definable compactness was introduced in \cite{steinhorn_definable_1996} as a property of definable spaces, however recent developments (namely \cite{andujar_guerrero_definable_2025}) treat it as a property of a more general class of spaces, that being definable topological spaces. 

Recall that, given structure $\cM$, a definable family of sets is a family of the form $\{\phi(\cM, b) : b \in \psi(\cM)\}$, for some formulas $\phi(\bar x,\bar y)$ and $\psi(\bar y)$, where $\phi(\cM, b)$ denotes the set of tuples $a \in \cM$ such that $\cM \models \phi(a,b)$ and, in an analogous manner, $\psi(\cM)$ is the set of tuples $a \in \cM$ such that $\cM \models \psi(a)$. The following was introduced in \cite{andujar_guerrero_definable_2025}.

\begin{definition*}
    Let $\cM$ be an o-minimal structure. A definable topological space $(X,\tau)$, with $X \subseteq M^n$ is a topological space such that there exists a definable family of subsets of $X$ which is a basis for $\tau$.
\end{definition*}

In particular, every definable space $X \subseteq M^n$ with definable charts is a definable topological space, as we may use the charts to transport the definable open boxes of $M^n$ to $X$.

\begin{definition*}
    Let $(X,\tau)$ be a definable topological space. We say that $X$ is definably compact if for every definable curve $\sigma:(a,b) \to X$, with $(a,b) \subseteq M$, there exists $\alpha,\beta \in X$ such that $\lim_{x \to a^+} \sigma(x) = \alpha$ and $\lim_{x \to b^-} \sigma(x) = \beta$.
\end{definition*}

\begin{remark}
    In \cite{andujar_guerrero_definable_2025}, this property is called ``\emph{curve-compactness}'', but in order to avoid confusion, we will simply name it ``\emph{definable compactness}''.
\end{remark}

 Given a definable set $X \subseteq M^n$, we denote by $S_X(M)$ the set of complete types $p \in S_n(M)$ such that $X \in p$, where here we are under common model-theoretic convention of identifying types with the
family of sets defined by formulas. Additionally, recall that a type 
$p$ is said to be definable if for every formula $\phi(x,y)$ there exists another formula 
$\psi(y)$ such that $\psi(\cM) = \{b \in M^{|y|}:\psi(x,b)\in p\}$. The following is part of Theorem A in \cite{andujar_guerrero_definable_2025}.

\begin{fact} \label{Facto andujar about compactness}
    Let $(X,\tau)$ be a Hausdorff definable topological space in an o-minimal structure $\cM$. The following are equivalent:
    \begin{itemize}
        \item $(X,\tau)$ is definably compact;
        \item Every downward directed definable family of non-empty closed sets has non-empty intersection (filter-compactness);
        \item Every definable type $p \in S_X(M)$ has a limit; i.e., there is a point in the intersection of every closed set in p (type-compactness).
    \end{itemize}
\end{fact}

Similarly to compactness, definable compactness in Hausdorff definable topological spaces enjoys the following familiar properties, which we will use freely throughout.

\begin{lemma}
Let $\cM$ be an o-minimal structure. Then the following hold:
\begin{enumerate}
    \item Let $(X,\tau)$ be a definably compact Hausdorff definable topological space and $Y \subseteq X$ is a definable subset. Then $Y$ is definably compact if and only if it is closed in $X$;
    \item Let $f:X \to Y$ is a continuous definable map between Hausdorff definable topological spaces. If $X$ is definably compact, then $f(X)$ is definably compact;
    \item  Let $X$ and $Y$ be two Hausdorff definable topological spaces. Then $X \times Y$ is definably compact if and only if both $X$ and $Y$ are definably compact.
\end{enumerate}
\end{lemma}

\begin{proof}
(1) Suppose first that $Y$ is closed in $X$. Let $\mathcal C$ be a downward directed
definable family of non-empty closed subsets of $Y$. Since $Y$ is closed in $X$,
the members of $\mathcal C$ are closed in $X$. Hence, by filter-compactness of
$X$, $\mathcal C$ has non-empty intersection. Thus $Y$ is filter-compact, and therefore definably compact.

Conversely, suppose that \(Y\) is definably compact. Let \(x\in \overline Y\) be an element in the topological closure of $Y$, let \(\mathcal B\) be a definable basis for \(\tau\), and consider the definable downward directed family 
\[
\mathcal F_x=\{Y\cap B : B \in \mathcal B \text{ and } x \in B\},
\]
By Proposition 4.3 in \cite{andujar_guerrero_definable_2025}, this family
extends to a definable type \(p\in S_Y(\mathcal M)\). Since \(Y\) is definably compact,
it is type-compact, so \(p\) has a limit \(y\in Y\).

We claim that \(x = y\) and prove it by contradiction. Assume that $x \neq y$, and note that since \(X\) is Hausdorff, we can find disjoint open
neighbourhoods \(U\) of \(y\) and \(V\) of \(x\). Take a basic open set \(B\) with
\(x\in B\subseteq V\). Since \(Y\cap B\in p\), and $Y\cap B \subseteq \overline{Y \cap B}$ it follows that $\overline{Y \cap B} \in p$, by completeness of $p$. The fact that \(y\) is a limit of \(p\) implies $y\in \overline{Y \cap B} \subseteq \overline{B}$. However, this contradicts the fact that \(U\cap B=\emptyset\). Hence \(x=y\) and therefore $x \in Y$, from which it follows that \(Y\) is closed.

(2) Assuming that $X$ is definably compact, this is straightforward using filter-compactness by taking the preimage under $f$ of a downward directed definable family of non-empty closed subsets of $f(X)$. 

(3) If \(X\times Y\) is definably compact, then \(X\) and \(Y\) are definably compact using both the fact that the projection maps onto $X$ and $Y$ are continuous definable maps and (2).
 
Conversely, suppose that \(X\) and \(Y\) are definably compact. Let $\gamma:(a,b)\to X\times Y$ be a definable curve. Then \(\pi_X\circ\gamma\) and \(\pi_Y\circ\gamma\) are definable~curves in \(X\) and \(Y\), and by definable compactness of \(X\) and \(Y\) both component curves have limits as \(t\to a^+\) and as \(t\to b^-\). By the definition of the product topology, \(\gamma\) has the corresponding product limits.
\end{proof}

\begin{remark}
The results proved in the remainder of the paper are stated for definable semigroups equipped with compatible Hausdorff definable space structures. However, several of the arguments only use the fact that the underlying topology is Hausdorff and definable in the sense of \cite{andujar_guerrero_definable_2025}, together with the basic properties of definable compactness for Hausdorff definable topological spaces recalled above. Consequently, many of the results below admit immediate formulations for Hausdorff definable topological semigroups.

We shall not pursue this additional level of generality here. The primary reason for this is that all natural examples of definable semigroups in o-minimal structures have a compatible definable space structure, and even a definable manifold structure (see Question \ref{question definable manifold}). Additionally, passing from the present setting to that of Hausdorff definable topological spaces would only require minor changes in terminology.
\end{remark}

\subsection{Preliminaries on Semigroups}
Finally, we will briefly recall a few definitions and results from semigroup theory that will be used throughout the paper. This is not intended to be a introduction to the subject; see \cite{howie_fundamentals_1996} for a classical reference.

Recall that a semigroup is a pair $(S,\cdot)$ where $S$ is a set and $\cdot:S \times S \to S$ is an associative operation on $S$. In a semigroup $S$, an element $e \in S$ is called an idempotent if $e^2 = e$, and the set of idempotents of $S$ is denoted by $E(S)$.

Given a semigroup $S$, we say that a non-empty subset $I$ is a left ideal of $S$ if, for all $i \in I$ and $s \in S$, we have $si \in I$. Right ideals are defined analogously and a subset is an ideal if it is a right and left ideal.

A monoid $M$ is a semigroup with an identity, i.e. there exists an element $e \in M$ such that for all $m \in M$, $em = me = m$.
Given a semigroup $S$, we define its induced monoid, denoted by $S^1$ as: If $S$ is already a monoid, then $S^1 = S$. If $S$ is not a monoid, then consider the set $S\cup \{1\}$, where $1$ is a new element and define the product $*$ in this set as: 
$$s * t = \begin{cases}
    s \cdot t & \text{if $s,t \in S$} \\
    s & \text{if $t = 1$} \\
    t & \text{if $s = 1$}
\end{cases}.$$

With this notation in mind, the left ideal generated by an element $s \in S$ is the subset $\es{s}$, the right ideal is $\di{s}$ and the ideal generated by $s$ is $\bi{s}$. 

There are five equivalence relations we can define on any semigroup, known as Green's equivalence relations, which play a key role in the study of the algebraic properties of semigroups. For each $s,t \in S$ define:
    \begin{align*}
        s \mathrel{\Lg} t\hspace{10pt} & \text{if}\hspace{10pt} \es{s}= \es{t};\\
        s \mathrel{\Rg} t\hspace{10pt} & \text{if}\hspace{10pt} \di{s}= \di{t};\\
        s \mathrel{\Jg} t\hspace{10pt} & \text{if}\hspace{10pt} \bi{s}= \bi{t}.
    \end{align*}
In other words, two elements are $\Lg$-related if they generate the same left ideal, $\Rg$-related if they generate the same right ideal and $\Jg$-related if they generate the same ideal.

Note that these are definable equivalence relations on $S$. 
Furthermore, we define the fourth Green's equivalence relation $\Hg$ as $\Hg = \Lg \cap \Rg$.

Recall that given two equivalence relations $\sigma,\rho$ on a set $X$, we define their composition $\sigma\rho$ as the equivalence relation given by $$(x,y) \in \sigma\rho \Leftrightarrow \exists z \in X : (x,z) \in \sigma \wedge (z,y) \in \rho.$$

We define the fifth and final Green's equivalence relation  $\Dg$ as  $\Dg := \Rg\Lg = \Lg\Rg$, i.e. given $s,t \in S$, we have $s \mathrel{\Dg} t$ if and only if there is $c \in S$ such that $s \mathrel{\Lg} c \mathrel{\Rg} t$. (For the proof that $\Rg\Lg = \Lg\Rg$ see Proposition 2.1.3 in \cite{howie_fundamentals_1996}).

The relationship between these five equivalence relations is given by the following Hasse diagram:

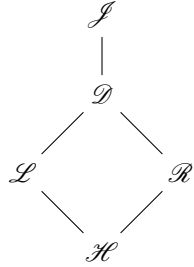
\begin{center}
    \begin{tikzpicture}[scale = .7]
    \node (J) at (0,3) {$\Jg$};
    \node (D) at (0,1.5) {$\Dg$};
    \node (L) at (-1.5, 0) {$\Lg$};
    \node (R) at (1.5, 0) {$\Rg$};
    \node (H) at (0, -1.5) {$\Hg$};
    \draw (J) -- (D) -- (L) -- (H) -- (R) -- (D);
\end{tikzpicture}
\end{center}

For a given $s \in S$, it is usual to denote by $L_s,R_s,J_s,H_s$ and $D_s$ the $\Lg,\Rg,\Jg,\Hg$ and $\Dg$ classes of $s$ respectively.

The following is an important result about $\Hg$-classes that will be essential later on. For more details, see Theorem 2.16 in \cite{howie_fundamentals_1996}.
\begin{fact}\label{H class with idempotent is group}
    In a semigroup $S$, an $\Hg$-class is a group if and only if it contains an idempotent.
\end{fact}

By a simple semigroup, we mean a semigroup with no proper ideals. Given a semigroup $S$ and $e,f \in E(S)$, we write $e\leq f$ if $ef = fe = e$. In the literature, this order is usually called the natural partial order on $E(S)$. Furthermore, we say that an idempotent $e$ is primitive if for all $f \in E(S), f\leq e \Rightarrow f = e$.

\begin{definition*}
\index{completely simple semigroup}
    We say that a semigroup $S$ is completely simple if it is simple and has at least one primitive idempotent.
\end{definition*}
 
\begin{example}\label{exemplo completamente simples}$ $
\vspace{-4pt}
    \begin{itemize}
            \item Any group is completely simple;
            \item  A semigroup $S$ is called a \emph{band} if $E(S) = S$. A band $S$ is called a \emph{rectangular band} if for all $a,b \in S$, we have $aba = a$. Any rectangular band is completely simple.
            \item If $G$ is a group and $B$ is a rectangular band, $G \times B$ is completely simple. Such semigroups are called \emph{rectangular groups}.
        \end{itemize}
\end{example}
An important class of completely simple semigroups are Rees matrix semigroups.

\begin{definition*} [Rees Matrix Semigroup] $ $
\index{$\mathscr{M}(I,S,\Lambda, P)$}
    Let $S$ be a semigroup, $I,\Lambda$ be non-empty sets, and let $P:\Lambda \times I \to S$ be a function.

    The Rees matrix semigroup $\mathscr{M}(I,S,\Lambda, P)$ is defined as the set $I \times S \times \Lambda$ with operation given by $$(i,s,\lambda) \cdot (j,t,\mu) = (i,s P(\lambda, j)t, \mu).$$
\end{definition*}

When dealing with Rees matrix semigroups, it is usual to denote $P(\lambda,i)$ by $p_{\lambda i}$.

The importance of Rees matrix semigroups comes from the fact that they are the only completely simple semigroups, up to isomorphism. For more details see Theorem 3.3.1 in \cite{howie_fundamentals_1996}.

\begin{theorem*}[Rees's Theorem] 
    Let $S$ be a semigroup. $S$ is completely simple if and only if there exists a group $G$, non-empty sets $I,\Lambda$ and a map $P:\Lambda \times I \to G$ such that $S\simeq \M I G \Lambda P$.
\end{theorem*}

\begin{example}
    Following the examples given in Example \ref{exemplo completamente simples}, we have that:

    \begin{itemize}
        \item Given any group $G$, then $G \simeq \M I G \Lambda P$, where $I = \Lambda = \{\bullet\}$ are sets with just one element, and $P:\Lambda \times I \to G$ maps $(\bullet,\bullet)\mapsto e$;

        \item Given a rectangular band $B$, it is a well known fact in semigroup theory that there are sets $I$ and $\Lambda$ such that $B \simeq I \times \Lambda$ where the product on $I \times \Lambda$ is given by $$(i,\lambda)(j,\mu) = (i,\mu).$$
        If we take $G = \{e\}$ to be the trivial group, then $B \simeq \M I G \Lambda P$ where $P$ is the only map from $\Lambda \times I$ to $G$;

        \item Given a rectangular group $G \times B$, then $G \times B \simeq \M I G \Lambda P$, where $B \simeq I \times \Lambda$ as described in the previous item, and $P$ maps any element of $\Lambda \times I$ to $e$.
    \end{itemize}
\end{example}

We can think of a completely simple semigroup $S \simeq \M I {G} {\Lambda} P$ as a grid

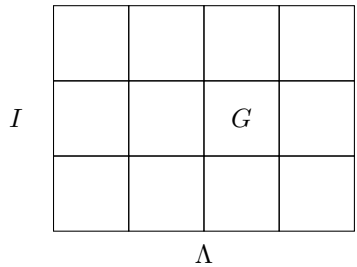
\tikzset{square matrix/.style={
    matrix of nodes,
    column sep=-\pgflinewidth, row sep=-\pgflinewidth,
    nodes={
      draw,
      minimum height=#1,
      anchor=center,
      text width=#1,
      align=center,
      inner sep=0pt
    },
  },
  square matrix/.default=1cm
}

\begin{center}
    \begin{tikzpicture}
\matrix[square matrix, nodes in empty cells,ampersand replacement=\&]
{
 \&  \& \& \\
 \&  \& $G$ \&  \\ 
 \& \& \& \\ 
};
\node at (0,-1.8) {{$\Lambda$}};
\node at (-2.5,0) {{$I$}};
\end{tikzpicture}

\end{center}

where each square is an isomorphic copy of $G$, the rows of the grid are indexed by the elements of $I$ and the columns by elements of $\Lambda$. As such, an element $(i_0,g_0,\lambda_0) \in \M I {G} {\Lambda} P$ corresponds to a choice of square in the grid, namely the one with coordinates $(i_0,\lambda_0)$ and then $g \in G$ is an element in the copy of $G$ responding to the square that we have picked.

Finally, the following well known fact
follows easily from Rees’s Theorem.
\begin{fact}\label{every idempotent primitive}
    In a completely simple semigroup, every idempotent is primitive.
\end{fact}

\section{Existence of idempotents and minimal ideals}

As said in the Section 1, from here onwards, we fix the following assumption:

\begin{quote}
    For the rest of the paper, unless stated otherwise, we will assume that every definable semigroup [group] we talk about and work with is equipped with a Hausdorff definable space structure that is compatible with the operations in said semigroup [group], i.e. the topology arising from such definable space structure makes it a topological semigroup [group].
\end{quote}

A classic result in the theory of topological semigroups is that any compact topological semigroup has at least one idempotent (see \cite{carruth_theory_1983}). As we will see now, definable compactness is enough to ensure the existence of idempotents.

\begin{theorem}\label{def compact semigroup has idempotent}
    Let $S$ be a definably compact definable semigroup. Then $S$ has at least one  idempotent.
\end{theorem}
\begin{cproof}
    This proof is done by considering three different cases.\\

    $\textbf{Case 1:}$ $S$ is commutative.

      Consider the definable family $\{sS : s \in S\}$ closed subsets parametrized by $S$.  Note that, as $S$ is commutative this family  is downward directed, i.e. for any $a,b \in S$, there exists $c \in S$ such that $aS \supseteq cS$ and $bS \supseteq cS$ (indeed, we may take $c=ab$). As such, by filter-compactness (Fact \ref{Facto andujar about compactness}), the intersection $\bigcap_{s \in S} sS$ is non-empty. 
    
    Let $T = \bigcap_{s \in S} sS$ and fix $t \in T$ and note that, as $T$ is a semigroup, $tT \subseteq T$. 
    
    \begin{claim*}
        The equality $T = tT$ holds.
    \end{claim*}

    We start by noting that $t\bigcap_{s \in S}sS = \bigcap_{s \in S}tsS$. Indeed, the inclusion $t\bigcap_{s \in S}sS \subseteq \bigcap_{s \in S}tsS$ follows easily. On the other hand, let $z  \in \bigcap_{s \in S}tsS$. Then, in particular, by taking $s$ to be $t$ in this intersection, we have that $z \in ttS$ and as such, there exists some $b \in S$ such that $z = ttb$. Write $a := tb$ and note that because $t \in \bigcap_{s \in S} sS$ (as $t \in T$), then $$a = tb \in \left(\bigcap _{s \in S} sS\right) b \subseteq \bigcap _{s \in S} sS.$$ So, $a \in \bigcap_{s \in S} sS$ and from the fact that $z = ta$ it follows that $z \in t \bigcap_{s \in S} sS$.
    
   As such,
    $$tT = t\bigcap_{s \in S}sS = \bigcap_{s \in S}tsS \supseteq T,$$
    meaning that $T \subseteq tT$ and therefore $tT = T$.
    
    \vspace{10pt}
    
    From this, we conclude that there exists $e \in T$ such that $te = t$ and $s \in T$ such that $ts= e$. So $e^2 = ets  = tes = ts = e$.\\

    $\textbf{Case 2:}$ $Z(S) \neq \emptyset$, where $Z(S) = \{x \in S : \forall y \in S, xy = yx\}$ is its center.

    Note that $Z(S)$ is a commutative definable subsemigroup of $S$. Indeed, if $x,y \in Z(S)$, then for any $s \in S$, we have $xys = xsy = sxy$, so $xy \in Z(S)$.

    Furthermore, note that $Z(S)$ is closed in $S$. To see this, let $C(a)$ denote the centralizer of an element $a \in S$, i.e. $C(a) = \{x \in S: xa = ax\}$. Then $Z(S) = \bigcap_{a \in S} C(a)$ and as such, it suffices to verify that each $C(a)$ is closed. For this, consider the definable continuous map $f_a:S \to S\times S$ given by $f(x) = (ax,xa)$. By definition, $C(a) = f^{-1}_a(\Delta)$, where $\Delta := \{(s,s) : s \in S\}$ is the diagonal in $S \times S$. As $S$ is Hausdorff, $\Delta$ is closed and thus $C(a)$ is a closed.

    As such, $Z(S)$ is a commutative definably compact definable semigroup, and thus it has an idempotent by Case 1.\\

    $\textbf{Case 3}$: General case.

    Fix $a \in S$, and consider its centralizer $C(a)$. Note that this is a definable semigroup, as if $x,y \in C(a)$, then $xya = xay = axy$. In particular, as we proved in Case 2, it is closed, from which it follows that $C(a)$ is a definably compact definable semigroup. 
    
    Finally, note that, by definition, $Z(C(a)) \neq \emptyset$ as $a \in Z(C(a))$, which implies that $C(a)$ has at least one idempotent by case 2.
\end{cproof}

Following the notation used in case 3 of this proof, given a set $A$, $\Delta_A$ from now on will always denote the diagonal on $A \times A$. When the set $A$ is implicit from the context, we will simply denote the diagonal by $\Delta$.

Let $(S,\cdot)$ be a topological semigroup. If $S$ is also a group, i.e. $S$ has inverses with respect to the product $\cdot$, it might be the case that the map $x \mapsto x^{-1}$ is not continuous, even though the product is. However, it is a well known fact that this is not possible if the topology in $S$ is compact, i.e. any compact topological semigroup $S$ that is algebraically a group, is a topological group (see \cite{carruth_theory_1983}). As we will see now, definably compact semigroups exhibit the same behaviour. 

\begin{lemma}\label{semigroup def compact + group = topological group}
    Let $S$ be a definably compact definable semigroup. If $S$ is algebraically a group, then it is a definably compact definable group, i.e. the inverse map is also continuous.
\end{lemma}
\begin{cproof}
    Let $m:S \times S \to S$ denote the multiplication in $S$, $i: S \to S$ the inverse function and $e$ its identity.

    Start by noting that the graph of the inverse function $\Gamma(i) = \{(x,x^{-1}) : x \in S\}$ is given by: $$\Gamma(i) = \{(x,y) \in S \times S: m(x,y) = e\} = m^{-1}(\{e\}).$$

    By continuity of $m$, $\Gamma(i)$ is closed.

    Let $C\subseteq S$ be a definable closed subset and let $\pi_1: S \times S \to S$ denote the projection into the first coordinate. Note that 
    $$i^{-1}(C) = \pi_1((S \times C) \cap \Gamma(i))$$,

    and additionally that $(S\times C) \cap \Gamma(i)$ is closed in $S\times S$, hence definably compact. As the image of a definably compact set under a definable continuous map is definably compact, $\pi_1((S \times C) \cap \Gamma(i)) = i^{-1}(C)$ is definably compact and therefore closed.

    As $S$ is equipped with a definable space structure, any closed set in $S$ is the intersection of definable closed subsets of $S$. Hence, it follows that the pre-image of any closed set under $i$ is closed, meaning that $i$ is continuous.
\end{cproof}

\begin{remark}
    If $S$ is equipped with a definably compact definable manifold structure, and $S$ is algebraically a group, by Lemma \ref{semigroup def compact + group = topological group}, $S$ will be a topological group. This implies that $S$ is equipped with Pillay's topology, as it is the unique manifold topology on any definable group that makes it a topological group.
\end{remark}

The existence of idempotents provided by Theorem \ref{def compact semigroup has idempotent} combined with Lemma \ref{semigroup def compact + group = topological group} allows us to prove the following results.

\begin{proposition} \label{def compact monoid com 1 idempotente é grupo}
    Let $S$ be a definably compact monoid such that $E(S) = \{1\}$, where $1$ is the identity of $S$. Then $S$ is a definably compact definable group.
\end{proposition}
\begin{cproof}
    Let $s \in S$ and note that $sS$ is a definably compact definable semigroup. By Theorem \ref{def compact semigroup has idempotent}, $sS \cap E(S) \neq \emptyset$ , i.e. $1 \in sS$. This implies $s$ has a right inverse. Analogously, we can prove that $s$ has a left inverse. 

    Moreover, the right and left inverses of $s$ have to be equal. To see this, let $l$ be the left inverse of $s$ and $r$ the right inverse. Then $$l = l1 = lsr = 1r = r,$$

    meaning that, algebraically, $S$ is a group.
    
    The inverse operation is obviously definable, meaning that $S$ is a definable group. By Lemma \ref{semigroup def compact + group = topological group}, $S$ is a topological group with respect to its topology, and thus a definably compact definable group.
\end{cproof}

Note that in this proposition, definable compactness really is necessary to conclude that $S$ is a group. For example, let $\cM = (\R, +, \cdot, 0, 1)$ and consider the definable semigroup $S = \{x \in \R : x \geq 1\}$ with the operation being given by the product. Then, $S$ is a monoid whose only idempotent is $1$, however, it is not a group.

\begin{lemma} \label{def compact subsemigroup of def group is subgroup}
    Let $G$ be a definable group. Then every definably compact subsemigroup of $G$ is a subgroup.
\end{lemma}
\begin{cproof}
    Let $S\subseteq G$ be a definably compact subsemigroup of $G$. By Theorem \ref{def compact semigroup has idempotent}, $S$ has an idempotent, and as the only idempotent in $G$ is the identity, we have that $E(S) = \{e\}$. As $e \in S$, $S$ is actually a monoid. By Proposition \ref{def compact monoid com 1 idempotente é grupo}, $S$ is a subgroup.
\end{cproof}

With this lemma, we can prove the following result, which is direct generalization of Theorem 5.1 in \cite{peterzil_generic_2007}.
\begin{proposition}\label{def subsemigroups of def compact are groups}
    Let $G$ be a definable group, $S$ be a definable subsemigroup of $G$, and let $\bar S$ denote the closure of $S$. If $\bar S$ is definably compact, then $S$ is a subgroup of $G$ and  $S = \bar S$.

    Additionally, if $G$ itself is definably compact, then every definable subsemigroup of $G$ is a subgroup.
\end{proposition}
\begin{cproof}
    Let $G$ be a definable group, $S$ be a definable subsemigroup such that $\bar S$ is definably compact. Start by noting that $\bar S$ is a subsemigroup of $G$, and thus, by Lemma \ref{def compact subsemigroup of def group is subgroup}, we get that $\bar S$ is a subgroup.  By Theorem 1.8 of Chapter 4 in \cite{dries_tame_1998}, we know that $\dim(\bar S \setminus S) < \dim S$, i.e. $S$ is large in $\bar S$. Let $x \in \bar S$, and note that by Lemma 2.1 of \cite{pillay_groups_1988}, there exist generic elements $a$ and $b$ of $\bar S$ such that $x = a\cdot b$. If either $a$ or $b$ where in $\bar S \setminus S$, then we would have that $\dim(\bar S \setminus S) = \dim S$, which is a contradiction. As such, $a,b \in S$, meaning that:
    $$\bar S \subseteq S \cdot S \subseteq S \subseteq \bar S.$$
    Therefore, $S = \bar S$.

    Assume now that $G$ is itself definably compact. As $\bar S$ is closed, it is definably compact and we can apply the first half of this corollary to conclude that $S$ is a subgroup.
\end{cproof}

Recall that a semigroup $S$ is said to be right cancellative if for all $a,b,c \in S$, $ac = bc$ implies that $a = b$. Left cancellative semigroups are defined analogously and a cancellative semigroup is simply a right and left cancellative semigroup. 

It is very straightforward to check that any group is necessary cancellative, but the converse is not true.

\begin{example}
        Fix the structure $\mathcal M =(\mathbb R, +, -, \cdot, 0 , 1)$. Consider the definable semigroup $S = (\mathbb R_{>0}, +)$. Then $S$ is cancellative, but it isn't a group (in fact, it isn't even a monoid).
\end{example}

\begin{corollary}\label{def compact semigroup is group iff it is cancelative}
    Let $S$ be a definably compact semigroup. Then $S$ is a group if and only if it is cancellative.
\end{corollary}
\begin{cproof}
    By Theorem \ref{def compact semigroup has idempotent}, $S$ has some idempotent $e$, which we fix. For any $x \in S$, we have that $ex = e^2x = e(ex)$ and so $x = ex$. Analogously, $xe = xe^2 = (xe)e$ which implies that $x = xe$.

    Therefore, $e$ is an identity of $S$ and $E(S) = \{e\}$, making $S$ a definably compact monoid with only one idempotent. By Proposition \ref{def compact monoid com 1 idempotente é grupo}, $S$ is a definable group.
\end{cproof}

Corollary \ref{def compact semigroup is group iff it is cancelative} thus shows that cancellativity is the property that `picks' the groups among the definably compact semigroups.

Let $S$ be a semigroup. By a minimal ideal, we mean an ideal $I$ such that, for any other ideal $J$, if $J \subseteq I$ then $J = I$. Equivalently, an ideal $I$ is minimal if for any other ideal $J$, we have that $I \subseteq J$. A left ideal $I$ is said to be a minimal left ideal if it is minimal among the left ideals, and a minimal right ideal if it is minimal among the right ideals.

A key property of compact semigroups is that such semigroups always have minimal left and right ideals, and a unique minimal ideal (see \cite{abodayeh_compact_1997}). As we will see now, this is also the case for definably compact definable semigroups.

We start by proving the existence of a unique minimal ideal.

\begin{theorem}\label{unique minimal ideal}
    Let $S$ be a definably compact definable semigroup. Then $S$ has a unique minimal ideal. Furthermore, that minimal ideal is definable and definably compact.
\end{theorem}
\begin{cproof}
    For each $a \in S$, let $J(a)$ denote the ideal $SaS$, which is closed. Then $\{J(a) : a \in S\}$ is a definable family of closed subsets parametrized by $S$. Furthermore, given $a_1,\ldots ,a_n \in S$, we have that $J(a_1\ldots a_n)\subseteq J(a_1) \cap \ldots  \cap J(a_n)$, which means that the family $\{J(a) : a \in S\}$ is downward directed. By filter-compactness (see Fact \ref{Facto andujar about compactness}) the intersection $I =\bigcap\{J(a) : a \in S\}$ is non-empty, and thus an ideal.

    \begin{claim*}
        $I$ is a minimal ideal.
    \end{claim*}
    To see this, let $J\subseteq I$ be an ideal, and let $x \in J$. Then $SxS\subseteq J \subseteq I \subseteq SxS$, meaning that $J = I$.

    For the uniqueness, let $I,J$ be two minimal ideals. As $IJ\subseteq I \cap J$ we get that $I = IJ = J$.

    Finally, note that $I = \bigcap_{a \in S} S a S$ is definable: indeed, it is defined by the formula $$\phi(x) := \forall a \in S \ \exists v,w \in S  (x = vaw).$$
    Additionally, because it is the intersection of closed subsets, it is closed and therefore definably compact.
\end{cproof}

In the literature, the unique minimal ideal of a compact semigroup $S$ is usually named the kernel of $S$, and it is usually denoted by $K(S)$ or simply $K$ when the semigroup $S$ is implicit from context.

We now establish the existence of minimal left and right ideals, starting with some technical lemmas.

The following is a generalization of Wallace's Swelling Lemma (see Theorem 1.9 in \cite{carruth_theory_1983}).

\begin{lemma}[Swelling Lemma] \label{swelling lemma}
    Let $S$ be a definably compact definable semigroup, and $A$ be a definable closed subset of $S$. Given $t \in A$, if $A\subseteq tA$, then $A = tA$.
\end{lemma}
\begin{cproof}
    Fix $t \in A$ with $A\subseteq tA$, and let $T$ denote the set $\{x \in S: tA\subseteq xA\}$. If $x,y \in T$, then $tA\subseteq xA\subseteq xtA\subseteq xyA$, meaning that $xy \in T$, i.e. $T$ is a definable subsemigroup of $S$.
    
    \begin{claim*}
        $T$ is closed.
    \end{claim*}
    Note that it is enough to show that $T$ is definably compact. Let $\gamma:(c,d) \to T$ be a definable curve, with $c,d \in M$. As $S$ is definably compact, there exists $\beta \in S$ such that $\lim_{x \to d^-} = \beta$. We claim that $\beta \in T$, i.e. $tA \subseteq \beta A$. Fix $ta \in tA$ and start by noting that, for each $x \in (c,d)$, $tA\subseteq \gamma(x)A$, meaning that $ta \in \gamma(x)A$. For each $x \in (c,d)$, let $A_x := \{y \in A : ta = \gamma(x)y\}$, which is uniformly defined by the formula $ y \in A \wedge ta = \gamma(x)y$. Furthermore, note that each $A_x$ is closed, as it is the preimage of the singleton $\{ta\}$ under the definable continuous maps $y\mapsto \gamma(x)y$. By Lemma 5.10 in \cite{andujar_guerrero_definable_2025}. there exists a definable map $\sigma:(c,d) \to S$ such that $\sigma(x) \in  A_x$ for all $x \in (c,d)$. 
    Because $S$ is definably compact, there exists $\alpha \in S$ such that $\lim_{x \to d^-} \sigma(x) = \alpha$. We claim that $ta = \beta \alpha$. For all $x \in (c,d)$, $ta = \gamma(x)\sigma(x)$, meaning that $$ta = \lim_{x \to b^-}(\gamma(x)\sigma(x)) = \lim_{x \to b^-} \gamma(x)\lim_{x \to b^-}\sigma(x) = \beta \alpha.$$
    So $\beta \in T$ and therefore $T$ is closed.\\

    By Theorem \ref{def compact semigroup has idempotent}, $T$ has some idempotent $e$. By the definition of $T$, we have that $A\subseteq tA \subseteq eA$. Fix $a \in A$, then there exists $b \in A$ such that $a = eb$ and hence, $ea = e(eb) = eb = a$, meaning that $eA \subseteq A$. From the inclusion $A\subseteq tA \subseteq eA$ that we had previously, the equality $A = tA$ follows.
\end{cproof}

\begin{proposition} \label{prop ideais minimais}
    Let $S$ be a definably compact definable semigroup and fix $e \in E(S)$. Then the following are equivalent:
    \begin{enumerate}
        \item $Se$ is a minimal left ideal;
        \item $eSe$ is a group;
        \item $eS$ is a minimal right ideal;
        \item $K = SeS$.
    \end{enumerate}
\end{proposition}
\begin{cproof}
This is Theorem 1.23 in \cite{carruth_theory_1983} with a few modifications.

Items $(1),(2)$ and $(3)$ can be proven to be equivalent only using algebraic arguments, and require no topology at all. The same goes for the proof that $(3)$ implies $(4)$. In these cases, see the proof of Theorem 1.23 of \cite{carruth_theory_1983}.

As for the proof that $(4)$ implies $(3)$, let $R$ be a right ideal contained in $eS$. Fix $x \in R$ and note that $xS \subseteq R\subseteq eS$, meaning that $SxS \subseteq SeS$. By minimality, $SxS = SeS$, and therefore, $e = axb$ for some $a,b \in S$. Now note that $$xS \subseteq R \subseteq eS = axbS\subseteq axS.$$
As $xS$ is a definable closed subset of $S$, by the Swelling Lemma (\ref{swelling lemma}), $xS = axS$, and by the previous chain of inclusions, it follows that $R = eS$, thus proving the minimality of $eS$.
\end{cproof}

With this, we can conclude the following.

\begin{corollary}\label{existem ideais minimais}
    Let $S$ be a definably compact definable semigroup. Then $S$ has minimal right and left ideals.
    
    In fact, $E(K) \neq \emptyset$, and for each $e \in E(K)$ we have that:

    \begin{enumerate}
        \item $Se$ is a minimal left ideal;
        \item $eSe$ is a group;
        \item $eS$ is a minimal right ideal;
        \item $K = SeS$
    \end{enumerate}

    Furthermore, each minimal right and minimal left ideal is definable and definably compact.
\end{corollary}
\begin{cproof}
    By Theorem \ref{unique minimal ideal} $K$ is definably compact, and thus by Theorem \ref{def compact semigroup has idempotent}, $K$ has some idempotent.
    The rest of the items follow easily from Proposition \ref{prop ideais minimais}.

    Additionally, if $L$ is a minimal left ideal, given $x \in L$, then $Sx = L$ by minimality, meaning that $L$ is both definable and definably compact. The same follows for any right ideal. 
\end{cproof}

By Theorem 3.2 in \cite{clifford_semigroups_1948}, if a semigroup has at least one minimal right and minimal left ideal, then its kernel is completely simple. From this fact, we obtain the following.

\begin{corollary}
    Let $S$ be a definably compact definable semigroup. Then its kernel is completely simple.
\end{corollary}

For further algebraic properties of the kernel and minimal ideals, see \cite{abodayeh_compact_1997}.

Given a semigroup $S$ and an ideal $I$, the Rees quotient of $S$ over $I$, denoted by $S/I$ is the set of equivalence classes of $S$ modulo the congruence $I^2 \cup \{(s,s) : s \in S\}$, i.e. every element of $I$ is identified together, and each element of $S \setminus I$ is identified only with itself. For any $i \in I$, we denote the equivalence class of $i$ in $S/I$ by $0$. This is because, for any $a \in S$, $[a]\cdot 0 = 0\cdot [a] = 0$ in $S/I$. Furthermore, as for any $a \in S \setminus I$, the equivalence class $[a]$ is set $\{a\}$, it is usual to simply denote this equivalence class by $a$ to simplify the notation.

If $S$ is a topological semigroup and $I$ is a closed ideal, then the Rees quotient $S/I$ is a topological semigroup with the quotient topology. There is a very useful way of characterizing the open sets in $S/I$: If $\pi:S \to S/I$ denotes the quotient map that sends each element of $S$ to each equivalence class, then the open sets in $S/I$ are precisely:
$$\{\pi(U) : \text{ $U\subseteq S$ is open and either $U \cap I = \emptyset$ or $I\subseteq U$}\}.$$

For more details, see \cite{khosravi_topological_2012}.

\begin{proposition}\label{rees quotient is def compact}
    Let $S$ be a definably compact definable semigroup and $I$ a definable closed ideal. Then $S/I$ definable. Furthermore, topologically $S/I$ with the quotient topology is a definably compact Hausdorff definable space.
\end{proposition}
\begin{cproof} $ $ 

\textbf{Step 1}: $S/I$ is definable.
    
    Fix an element $x \in I$ which we will denote by $0$. Consider the definable set $S \setminus I\cup \{0\}$ and the formula $\phi(a,b,c)$ given by 
    \begin{align*}
        &(a = 0  \wedge c = 0)\  \vee\\ &(b = 0 \wedge c = 0)\ \vee\\ &(a \in S\setminus I \wedge b \in S\setminus I \wedge((ab \in S \setminus I \wedge c = ab) \vee (ab \in I \wedge c = 0)).
    \end{align*}

    Then $\phi(a,b,c)$ defines a product $*$ in the following way: If either $a$ or $b$ is $0$, then $ab = 0$. On the other hand, if both $a$ and $b$ are in $S\setminus I$, then 
    $$a*b =\begin{cases}
        ab & \text{if $ab \in S \setminus I$}\\
        0 & \text{if $ab \in I$}
    \end{cases}. $$
    
    In particular, $(S\setminus I \cup \{0\},*)$ is isomorphic to the Rees quotient $S/I$ \\

    $\textbf{Step 2}$: $S/I$ is definable space.
    
    Let $(U_i,\phi_i)_{i=1,\ldots ,n}$ be a definable atlas of $S$. For each $i = 1,\ldots ,n$ let $V_i = U_i \cap I^c$, where $I^c = S \setminus I$ denotes the complement of $I$,  and consider the family $(V_i, \phi_i |_{V_i})_{i=1,\ldots ,n}$. Without loss of generality assume that $V_i \neq \emptyset$ for all $i = 1,\ldots ,n$. Denoting $\pi(I)$ by $0$, we get that the family   $(\pi(V_i), \phi_i |_{V_i} \circ \pi^{-1})_{i=1,\ldots ,n}$ is a definable atlas in $(S/I)\setminus \{0\}$. Adding an additional chart $\psi:\{0\} \to \{a\}$ for some $a \in M$, we get that $S/I$ is a definable space.\\

    $\textbf{Step 3}$: $S/I$ is Hausdorff and definably compact.

    We will start by proving that $S/I$ is Hausdorff. Let $a,b \in S/I$ with $a \neq b$. As each non-zero element in $S/I$ is an equivalence class with only one element, we will treat $a$ and $b$ as both an element of $S$ and as their equivalence classes in $S/I$. If neither $a$ nor $b$ is zero, as $S$ is Hausdorff, there are open disjoint open sets $U',V'$ such that $a \in U'$ and $b \in V'$. Let $U = U' \cap I^c$ and $V = V' \cap I^c$ and note that these are disjoint open sets that separate $a$ and $b$ in $S$. As $I \cap U$ and $I \cap V$ are empty, we have that $\pi(U)$ and $\pi(V)$ are disjoint open sets in $S/I$ with $a \in \pi(U)$ and $b \in \pi(U)$. 

    On the other hand, let $a = 0$ and $b \neq 0$, i.e. $b \not \in I$. By Lemma 3.5 of \cite{edmundo_o-minimal_2001}, there are open disjoint sets $U$ and $V$ such that $I \subseteq U$ and $b \in V$, therefore $\pi(U)$ and $\pi(V)$ are disjoint open sets in $S/I$ that separate $a =0$ and $b$.

    As for the curve completion, let $\gamma:(a,b) \to S/I$ be a definable curve. Assume without loss of generality that $0 \not \in \gamma((a,b))$. This induces a unique curve $\sigma:(a,b) \to S$ such that $\pi \circ \sigma = \gamma$. By definable compactness, there exists $\beta \in S$ such that $\lim_{t \to b^-} \sigma(t) = \beta$. Then, $$\pi(\beta) = \pi(\lim_{t \to b^-} \sigma(t)) = \lim_{t \to b^-} \pi(\sigma(t)) = \lim_{t \to b^-} \gamma(t).$$
    The proof for the case $t \to a^+$ is analogous.
\end{cproof}

\begin{remark}\label{remark about S and S/K}

    This proposition suggests that we can `split' the study of definably compact semigroups into the study of two families of semigroups: On one hand, any definably compact definable semigroup $S$ has a kernel $K$, which is a definably compact completely simple semigroup. On the other hand, if we quotient out the kernel and look at $S/K$, we will have a definably compact definable semigroup with zero.
\end{remark}

For the remainder of this paper, we shall focus on the first case and study definably compact completely simple semigroups.

\section{Definably compact completely simple Semigroups}

As previously noted in Section 3, the kernel of any definably compact semigroup is completely simple. This suggests that a thorough understanding of definably compact completely simple semigroups may provide a structural lens for arbitrary definably compact semigroups. The objective of this section is thus to study such semigroups.

We start by generalizing Rees's Theorem to include definability.

\begin{lemma}[Definable Rees's Theorem]\label{definable Rees theorem}
    Let $S$ be a completely simple semigroup definable in $\cM$. Then, there exist definable sets $I,\Lambda$, a definable group $G$ and a definable map $P:\Lambda \times I \to G$ such that $S$ is definably isomorphic to $\mathscr{M}(I,G,\Lambda,P)$.
\end{lemma}
\begin{cproof}
    This can be proven by merely following the classical proof of Rees's Theorem and noting that the sets $I,\Lambda$, the group $G$ and the map $P:\Lambda \times I \to G$ are actually definable. For this, we will outline the proof of Rees's Theorem given in \cite{howie_fundamentals_1996} (see Theorem 3.2.3).

    The construction goes as follows: Pick $e$ to be any idempotent in $S$. Set $I = Se \cap E(S)$ and $\Lambda = eS \cap E(S)$, which are definable. The group $G$ is set to be the $\Hg$-class of $e$, i.e. $G = H_e$, which is also definable. Furthermore, the map $P:\Lambda \times I \to G$ is given by $P(\lambda,i) = \lambda i$ which, again, is definable.

    Finally, the isomorphism between $S$ and $\mathscr{M}(I,G,\Lambda,P)$ is given by
    \begin{align*}
    \M {I} {G} {\Lambda} P &\to S \\
    (i,g,\lambda) &\mapsto ig\lambda,
\end{align*}
    which is definable.
\end{cproof}

This result characterizes the algebraic structure of any definable completely simple semigroup up to definable isomorphism. However, we are also interested in considering topologies on such semigroups. 

It is straightforward to prove that given any two topological spaces $I,\Lambda$, a topological group $G$ and a continuous map $P:\Lambda \times I \to G$, then $\mathscr{M}(I,G,\Lambda,P)$ equipped with the product topology on $I \times G \times \Lambda$ is a topological semigroup. However, not all topologies on completely simple semigroups arise in this way. For instance, we could have a topology on $I \times G \times \Lambda$ that is not a product topology but still makes $\mathscr{M}(I,G,\Lambda,P)$ a topological semigroup. This prompts the following definition, which appears on \cite{banakh_rees-suschkewitsch_2009}.

\begin{definition} \label{top paragroup}
\index{topological!paragroup}
    A topological completely simple semigroup $S$ is said to be a topological paragroup if there are topological spaces $I, \Lambda$, a topological group $G$ and a continuous map $P:\Lambda \times I \to G$ such that $S$ is topologically isomorphic (i.e. there exists an isomorphism that is also a homeomorphism) to $\M I G \Lambda P$ equipped with the product topology on the base set $I \times G \times \Lambda$.
\end{definition}

\begin{proposition} \label{def compact comp simples is top paragroup}
    Let $S$ be a definably compact completely simple semigroup. Then $S$ is a topological paragroup. 
    
    Moreover, $S \simeq \M I G \Lambda P$, where $I$ and $\Lambda$ are definable sets with a definably compact Hausdorff definable space structure, $G$ a definable group with a definably compact Hausdorff definable space structure making it a topological group, $P$ a definable continuous map, and the topological isomorphism is definable.
\end{proposition}
\begin{cproof}
    Following the proof of Lemma \ref{definable Rees theorem}, fix an idempotent $e \in S$ and let $I = Se \cap E(S)$, $\Lambda = eS \cap E(S)$, $G = H_e = eSe$ and $P:\Lambda \times I \to G$ given by $P(\lambda,i) = \lambda i$.

    Consider the map 
    \begin{align*}
        f:S &\to S \times S \\ s &\mapsto(s,s^2)
    \end{align*}

    and note that $E(S) = f^{-1}(\Delta)$, i.e. $E(S)$ is closed.
    As $I$ is the intersection of two closed subsets, it is itself closed and thus definably compact. The same argument shows that $\Lambda$ is also definably compact.
    
    Additionally, as $G = eSe$, we conclude that $G$ is definably compact. Furthermore, with the topology induced by $S$, $G$ is a topological semigroup. By Lemma \ref{semigroup def compact + group = topological group}, $G$ is a topological group. 

    Recall that the isomorphism from $\M I G \Lambda P$ to $S$ is given by 
    \begin{align*}
        \phi: \M I G \Lambda P &\to S \\
        (i,g,\lambda) &\to ig\lambda
    \end{align*}
        with its inverse being
        \begin{align*}
            \psi:S &\to \M I G \Lambda P\\
            s & \mapsto (s(ses)^{-1}, ses, (ese)^{-1}s),
        \end{align*}
        where the inverses are being taken inside of $G$. By continuity of the product and of the inverse map, both $\phi$ and $\psi$ are continuous, as well as $P$.
\end{cproof}

    In particular, the kernel of any definably compact definable semigroup $S$ is a topological paragroup.

Recall that in a semigroup $S$, an element $s$ is said to be regular if there exists $t \in S$ such that $sts = s$. A semigroup is said to be regular if every element is regular. The following is Proposition 2.4.2 of \cite{howie_fundamentals_1996}.

\begin{fact} \label{regular subsemigroup L,R,H classes}
    Let $S$ be a semigroup and $T$ be a regular subsemigroup of $S$. Then
    \begin{align*}
        &\Lg^T = \Lg^S \cap (T \times T),\\
        &\Rg^T = \Rg^S \cap (T \times T),\\
        &\Hg^T = \Hg^S \cap (T \times T),
    \end{align*}

    where $\mathscr{K}^S$ and $\mathscr{K}^T$ denotes the green relation $\mathscr{K}$ in $S$ and in $T$ respectively, where $\mathscr{K} \in \{\Lg, \Rg, \Hg\}$.
\end{fact}

With this proposition, we can prove the following.

\begin{proposition} \label{subsemigroup is completely simple}
    Let $S \simeq \M I G \Lambda P$ be a definably compact completely simple semigroup and $T$ be a definable subsemigroup of $S$. Then $T$ is completely simple. 
\end{proposition}
\begin{cproof}
     Let $T$ be a definable subsemigroup of $S \simeq \M I G \Lambda P$. Start by noting that for each $\Hg$-class $H_e^S$ of $S$, with $e \in E(S)$, the intersection $T \cap H_e^S$ is either a subsemigroup of $H_e^S$ or the empty set. In particular, if $T \cap H_e^S$ is not empty, then by Proposition \ref{def subsemigroups of def compact are groups}, it is a subgroup of $H_e^S$ (recall that by Fact \ref{H class with idempotent is group}, $H_e^S$ is a group). This means that the elements of $E(T)$ are precisely the idempotents of $S$ where $T\cap H_e^S$ is a group.

     In particular, note that we can write $T$ as a union of groups $T = \bigcup_{e \in E(T)} (T \cap H_e^S)$, which implies that $T$ is a regular semigroup. By Fact \ref{regular subsemigroup L,R,H classes}, for each $e \in E(T)$, $H_e^T = T \cap H_e^S$ and, in particular, because $T = \bigcup_{e \in E(T)} (T \cap H_e^S) = \bigcup_{e \in E(T)} H_e^T$, we conclude that every $\Hg$-class of $T$ is a group.
    \begin{claim*}
        $T$ is simple.
    \end{claim*}
     Let $I \subseteq T$ be an ideal and note that if $I \cap H^T \ne \emptyset$ for some $\Hg$-class $H^T$ of $S$, then $H^T \subseteq I$.

     Let $x \in I$ and let $e_0 \in E(T)$ be an idempotent such that $x \in H^T_{e_0}$, and then $H^T_{e_0} \subseteq I$. Now, let $f \in E(T)$ be any other idempotent, and note that $fe_0f \in H^T_f$. As $I$ is an ideal, $fe_0f \in I$ and as such, we conclude that $H_f^T \subseteq I$ for each $f \in E(T)$. As $T = \bigcup_{f \in E(T)} H_f^T$, we get that $I = T$ and thus $T$ is simple.\\

     Now, as $S$ is completely simple, every idempotent of $S$ is primitive (Fact \ref{every idempotent primitive}), meaning that every idempotent of $T$ is primitive as well.
\end{cproof}

This proposition gives us a simple characterization of the subsemigroups of any definably compact completely simple semigroup, as we will see now.

\begin{corollary}\label{Forma dos subsemigroups nos completamente simples}
     Let $S \simeq \M I G \Lambda P$ be definably compact completely simple semigroup and $T$ be a definable subsemigroup of $S$. Then there are definable subsets $J \subseteq I$, $\Gamma \subseteq \Lambda$ and a definable subgroup $W \subseteq G$ such that $T \simeq \M J W \Gamma {P|_{\Gamma \times J}}$.
\end{corollary}
\begin{cproof}
    Let $e \in E(T)$. Recall that by the proof of Lemma \ref{definable Rees theorem}, we may take $I = Se \cap E(S)$, $\Lambda = eS \cap E(S)$, $G = H_e$ and $P:\Lambda \times I \to G$ to be the map given by $P(\lambda,i) = \lambda i$.

    By Proposition \ref{subsemigroup is completely simple}, $T$ is completely simple. As such, following the same construction, taking $J = Te \cap E(T)$, $\Gamma = eT \cap E(T)$, $W = H_e^T$ (where $H^T_e$ is the $\Hg$-class of $e$ in $T$) and $P:\Gamma \times J \to W$ to be the map given by $Q(\gamma,j) = \gamma j$, we have that $T \simeq \M J W \Gamma Q$. Furthermore, it is straightforward to check that $J \subseteq I$, $\Gamma \subseteq \Lambda$, $W \subseteq G$ and $Q = {P|_{\Gamma \times J}}$.
\end{cproof}

With this corollary, we can generalize the properties in Fact \ref{facts about definable groups} to completely simple semigroups. However, as we will see shortly, most of them only work when both $|I|$ and $|\Lambda|$ are finite and the semigroup we are working with is definably compact.

Just some terminology before stating the result: we say that a definable semigroup has the descending chain condition (DCC) if any descending sequence of definable subsemigroups $T_1 \supseteq T_2\supseteq \ldots $ eventually stabilizes. 

\begin{proposition}\label{completely simple DCC}
Let $S \simeq \M I G \Lambda P$ be a definably compact completely simple semigroup. Then $S$ has the DCC if and only if both $I$ and $\Lambda$ are finite.
\end{proposition}
\begin{cproof}
    Assume that $S$ has the DCC, and for the sake of contradiction, assume that either $I$ or $\Lambda$ are infinite. Without loss of generality, assume that $I$ is infinite, and fix an infinite subset $\{i_1,i_2,\ldots \} \subseteq I$. For each $n \geq 1$, let $I_n = I \setminus \{i_1,\ldots ,i_n\}$, which is definable. Then $$\M {I_1} G \Lambda {P|_{\Lambda\times I_1}} \supset \M {I_2} G \Lambda {P|_{\Lambda\times I_2}} \supset\ldots $$
    is an infinite descending chain of definable subsemigroups that never stabilizes, contradicting the fact that $S$ has the DCC.

    On the other hand, assume that both $I$ and $\Lambda$ are finite. Let $T_1 \supseteq T_2 \supseteq \ldots$ be a descending chain of definable subsemigroups. By Corollary \ref{Forma dos subsemigroups nos completamente simples}, for each $i \geq 1$, let $T \simeq \M{J_i}{W_i}{\Gamma_i}{P|_{\Gamma_i \times J_i}}$, where $J_i \subseteq I$, $\Gamma_i \subseteq \Lambda$ and $W_i$ is a definable subgroup of $G$.

    As $I$ is finite, the chain $J_1 \supseteq J_2 \supseteq \ldots$ eventually stabilizes to a subset $J \subseteq I$. Analogously, the chain $\Gamma_1 \supseteq \Gamma_2 \supseteq \ldots$ eventually stabilizes to a subset $\Gamma \subseteq \Lambda$. As $G$ has the DCC, the chain of definable subgroups $W_1 \supseteq W_2 \supseteq \ldots$ eventually stabilizes to a subgroup $W$. As such, the original chain $T_1 \supseteq T_2 \supseteq \ldots$ stabilizes to $\M J W \Gamma P_{\Gamma \times J}$.
\end{cproof}

If a completely simple semigroup is isomorphic to $\M I G \Lambda P$ with $I$ and $\Lambda$ being finite sets, then it will only have finitely many $\Hg$-classes, and for that reason, we say that such semigroups is $\Hg$-finite.

Note that in this proposition, definable compactness is necessary. Consider for example the structure $\cM = (\R, <, +, \cdot, 0, 1)$. Then $(\R,+)$ is a definable completely simple semigroup isomorphic to $\M I \R \Lambda P$, where $I = \Lambda = \{0\}$ and $P(0,0) = 0$. Then $(0,\infty) \supset (1,\infty) \supset \ldots  \supset (n,\infty) \supset\ldots $ is an example of a strictly descending chain of definable subsemigroups.

As stated earlier in Remark \ref{remark about S and S/K}, to study the structure of general definably compact semigroups $S$, we need to both understand the structure of definably compact completely simple semigroups (the kernel of $S$) and definably compact semigroups with zero (as is the case of $S/K$). The following result is an example of this.

\begin{corollary}\label{DCC equiv}
    Let $S$ be a definably compact definable semigroup. Then $S$ has the DCC if and only if $K$ is $\Hg$-finite and $S/K$ has the DCC.
\end{corollary}
\begin{cproof}
    Let $\pi:S\to S/K$ denote the projection map that maps every element of $S$ to its equivalence class in $S/K$. In particular, note that $\pi|_{S\setminus K}:S\setminus K \to (S/K)\setminus \{0\}$ is a bijection.

    Start by assuming that $S$ has the DCC. 
    
    If $K$ is not $\Hg$-finite then, by Proposition \ref{completely simple DCC}, there would exist a strictly descending chain of definable subsemigroups in $K$, which would contradict the fact that $S$ has de DCC. 
    
    As for the quotient $S/K$, let $H_1 \supseteq H_2 \supseteq \ldots$ be a descending chain of definable subsemigroups of $S/K$. We prove that $S/K$ has the DCC by considering two cases:

    $\textbf{Case 1}$: There exists $N \geq 1$ such that $0 \not \in H_N$.

    If this is the case, then for any $n \geq N$, we have that $0 \not\in H_n$, which means that $\pi|_{S\setminus K}(H_N) \supseteq \pi|_{S\setminus K}^{-1}(H_{N + 1}) \supseteq \ldots$ is a descending chain of definable subsemigroups of $S$, and thus must stabilize as $S$ has the DCC. This implies that $H_N \supseteq H_{N + 1} \supseteq \ldots$ also stabilizes.

    $\textbf{Case 2}$: For all $n \geq 1$, $0 \in H_n$.

    For each $n \geq 1$, consider the definable subsemigroup $H'_n = \pi|_{S\setminus K}^{-1}(H_n\setminus\{0\}) \cup K$ of $S$. As $S$ has the DCC, the chain $H'_1 \supseteq H'_2 \supseteq \ldots$ eventually stabilizes, which implies that the original chain $H_1 \supseteq H_2 \supseteq \ldots$ stabilizes as well.\\

    On the other hand, assume that $K$ is $\Hg$-finite and that $S/K$ has the DCC. Let $H_1 \supseteq H_2 \supseteq \ldots$ be a descending chain of definable subsemigroups in $S$. 
    
    For each $n \geq 1$, let $T_n = H_n \cap K$. As $K$ is $\Hg$-finite, by Proposition \ref{completely simple DCC}, there exists $N \geq 1$ such that $n \geq N$ implies that $T_n = T_N$.

    Additionally, for each $n \geq 1$, let $W_n = \pi(H_n)$. Then, as $S/K$ has the DCC, there exists $M \geq 1$ such that for all $m \geq M$, we have that $W_m= W_M$.  

    Note that for each $n \geq 1$, $H_n = T_n \cup \pi|_{S\setminus K}(W_n \setminus \{0\})$. Therefore, for $n,m \geq \max\{N,M\}$, we have that $H_n = H_{m}$, meaning that $S$ has the DCC.
\end{cproof}

\begin{remark}
    A special case of this is when $S$ is commutative. A completely simple semigroup with non-empty centre is necessarily a group, meaning that in this case, the kernel $K(S)$ is a group, i.e. it only has one idempotent. As such, by Corollary \ref{DCC equiv}, $S$ has the DCC if and only if $S/K$ has the DCC.

    Another instance of this is when a definably compact semigroup $S$ has only finitely many idempotents to begin with. Then Corollary \ref{DCC equiv} also implies that $S$ has the DCC if and only if $S/K$ has the DCC.
\end{remark}

Note that if $S \simeq \M I G \Lambda P$ carries a definable manifold structure and is $\Hg$-finite, then the topology induced in $G$ is also a definable manifold topology. This induced topology is precisely Pillay's topology, as it is the unique definable manifold topology on any definable group that makes it a topological group. As such, the properties given by Remark \ref{facts about definable groups} apply to $G$. With this fact in mind, we have the following:

\begin{proposition} \label{H finite semigroup equivalencia}
    Let $S \simeq \M I G \Lambda P$ be a definably compact completely simple semigroup with a compatible definable manifold structure. The following are equivalent:
    \begin{enumerate}
        \item Both $I$ and $\Lambda$ are finite;
        \item There are only finitely many definable subsemigroups $T$ where $\dim T $ is equal to $\dim S$;
        \item A subsemigroup $T$ of $S$ is open if and only if $\dim T = \dim S$;
        \item Any definable subsemigroup of $S$ is closed.
    \end{enumerate}
\end{proposition}
\begin{cproof}$ $
$(1) \Rightarrow (2)$.  Assume that both  $I$ and $\Lambda$ are finite. For any subsemigroup $T \simeq \M J W \Gamma {P|_{\Gamma \times J}}$, note that $\dim(T) = \dim(W)$. Fixing a definable subgroup $W \subseteq G$, there are at most $2^{|I|}2^{|\Lambda|}$ subsemigroups of $S$ of the form $T \simeq \M J W \Gamma {P|_{\Gamma \times J}}$. This means that, if we let \(\mathcal{S}\) denote the set of definable subsemigroups \(T\) of \(S\) such that \(\dim T=\dim S\), and \(\mathcal{G}\) denote the set of definable subgroups \(W\) of \(G\) such that \(\dim W=\dim G\), then $
|\mathcal{S}| \leq 2^{|I|}\cdot |\mathcal{G}| \cdot 2^{|\Lambda|}.$ By Fact \ref{facts about definable groups}, this last quantity is finite.

$(2)\Rightarrow(1)$. For the sake of contradiction, assume that $I$ is infinite. For each $x \in I$, let $I_x = I \setminus \{x\}$. In particular, note that $\dim I_x = \dim I$. This means that, for each $x\in I$, $\M{I_x}{ G}{\Lambda}{P|_{\Lambda \times I_x}}$ is a subsemigroup with the same dimension as $S$, thus contradicting (2). The proof that $\Lambda$ is finite is analogous.

$(1)\Rightarrow(3)$. Let $T$ be a definable subsemigroup and write $T = \M J W \Gamma {P|_{\Gamma \times J}}$ as per Corollary \ref{Forma dos subsemigroups nos completamente simples}. Note that $\dim T = \dim(J \times W \times \Gamma) = \dim(J) + \dim(W) + \dim(\Gamma)$. As both $J$ and $\Gamma$ are finite, we get that $\dim T = \dim W$. On the other hand, $\dim S = \dim(I \times G \times \Lambda) = \dim G$ for the same reason.

    Now, assume that $T$ is open. As the projection maps are open, we get that $W$ is open in $G$. By Lemma 2.11 of \cite{pillay_groups_1988}, we conclude that $\dim W = \dim G$, and thus $\dim T = \dim S$.
    On the other hand, if $\dim T = \dim S$, then $\dim W = \dim G$. By Lemma 2.11 of \cite{pillay_groups_1988}, $W$ is open in $G$. As both $J$ and $\Gamma$ are open, we get that $T = J \times W \times \Gamma$ is open.

$(3) \Rightarrow(1)$. We will show that $I$ is finite, as the proof that $\Lambda$ is finite is analogous.

\textbf{Claim:} Let $X$ be a definable subset of $I$. If $\dim X = \dim I$, then $X$ is open in $I$.

 Let $X \subseteq I$ be a definable subset with $\dim X = \dim I$. Consider the definable semigroup $\M X G \Lambda {P|_{\Lambda \times X}}$. As $\dim(X \times G \times \Lambda) = \dim(S)$, by our hypothesis, $X \times G \times \Lambda$ is open, and as the projection map is open, $X$ is open in $I$.

This claim is the property that will force the finiteness of $I$. For the sake of contradiction assume that $I$ is infinite, i.e. $\dim(I) = d > 0$. 

Let $\{(U_i,\phi_i)\}_{i=1,\ldots,k}$ be a definable atlas compatible with the definable space structure of $I$. As $I = \bigcup_{i=1}^k U_i$, there exists some $U_i$ with $\dim U_i = d$, which without loss of generality, we assume it to be $U_1$.

As $\phi_1$ is a definable bijection, $\dim(\phi_1(U_1)) = d$. Let $C$ be a cell of dimension $d$ with $C \subseteq \phi_1(U_1)$ and let $\pi$ be a projection onto $M^d$ such that $\pi|_C:C \to \pi(C)$ is an homeomorphism into a cell of dimension $d$. Let $B = (a_1,b_1) \times \ldots\times (a_d,b_d)$ be some open box contained in $\pi|_C(C)$ (see image below).

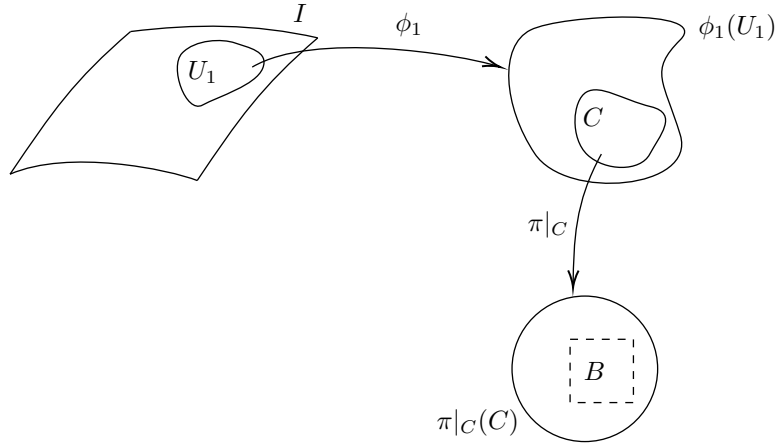
\begin{center}

\tikzset{every picture/.style={line width=0.5}} 

\begin{tikzpicture}[x=0.75pt,y=0.75pt,yscale=-.8,xscale=.8]

\draw    (80,115.9) .. controls (107.04,75.92) and (120,56.9) .. (156,25.9) ;
\draw    (156,25.9) .. controls (197,20.9) and (234,20.9) .. (274,29.9) ;
\draw    (198,119.9) .. controls (223.84,81.52) and (238,60.9) .. (274,29.9) ;
\draw    (80,115.9) .. controls (117,99.9) and (180,112.9) .. (198,119.9) ;
\draw   (188.48,39.52) .. controls (196.88,30.72) and (213.84,30.32) .. (221.84,33.12) .. controls (229.84,35.92) and (243.68,38.32) .. (239.28,49.12) .. controls (234.88,59.92) and (211.68,66.72) .. (201.28,72.32) .. controls (190.88,77.92) and (180.08,48.32) .. (188.48,39.52) -- cycle ;
\draw   (408.13,24.93) .. controls (428.13,14.93) and (522.53,11.13) .. (502.53,31.13) .. controls (482.53,51.13) and (492.13,56.93) .. (502.53,91.13) .. controls (512.93,125.33) and (430.8,132.93) .. (410.8,102.93) .. controls (390.8,72.93) and (388.13,34.93) .. (408.13,24.93) -- cycle ;
\draw   (443.33,64.13) .. controls (452.21,58.93) and (466.05,68.93) .. (484.05,73.33) .. controls (502.05,77.73) and (489.25,92.13) .. (483.25,103.33) .. controls (477.25,114.53) and (455.25,114.13) .. (444.05,103.73) .. controls (432.85,93.33) and (434.45,69.33) .. (443.33,64.13) -- cycle ;
\draw[dashed]   (433.07,220.07) -- (473,220.07) -- (473,260) -- (433.07,260) -- cycle ;
\draw   (396.47,238.7) .. controls (396.47,213.35) and (417.02,192.8) .. (442.37,192.8) .. controls (467.72,192.8) and (488.27,213.35) .. (488.27,238.7) .. controls (488.27,264.05) and (467.72,284.6) .. (442.37,284.6) .. controls (417.02,284.6) and (396.47,264.05) .. (396.47,238.7) -- cycle ;
\draw    (232.67,48.33) .. controls (271.28,25.52) and (362.87,40.16) .. (387.2,47.7) ;
\draw [shift={(388.93,48.27)}, rotate = 198.97] [color={rgb, 255:red, 0; green, 0; blue, 0 }  ][line width=0.75]    (10.93,-3.29) .. controls (6.95,-1.4) and (3.31,-0.3) .. (0,0) .. controls (3.31,0.3) and (6.95,1.4) .. (10.93,3.29)   ;
\draw    (452.8,103.2) .. controls (435.81,133.91) and (435.47,158.21) .. (434.96,185.9) ;
\draw [shift={(434.93,187.6)}, rotate = 271.08] [color={rgb, 255:red, 0; green, 0; blue, 0 }  ][line width=0.75]    (10.93,-3.29) .. controls (6.95,-1.4) and (3.31,-0.3) .. (0,0) .. controls (3.31,0.3) and (6.95,1.4) .. (10.93,3.29)   ;

\draw (256.8,6.8) node [anchor=north west][inner sep=0.75pt]    {$I$};
\draw (190.48,42.92) node [anchor=north west][inner sep=0.75pt]    {$U_{1}$};
\draw (512.4,13.4) node [anchor=north west][inner sep=0.75pt]    {$\phi _{1}( U_{1})$};
\draw (439.41,72.33) node [anchor=north west][inner sep=0.75pt]    {$C$};
\draw (347.33,262.4) node [anchor=north west][inner sep=0.75pt]    {$\pi |_{C}( C)$};
\draw (440,232.73) node [anchor=north west][inner sep=0.75pt]    {$B$};
\draw (321.33,12.4) node [anchor=north west][inner sep=0.75pt]    {$\phi _{1}$};
\draw (404.67,137.73) node [anchor=north west][inner sep=0.75pt]    {$\pi |_{C}$};

\end{tikzpicture}
\end{center}

Then, $\dim \pi|_C(C) \setminus B = \dim\pi|_C (C)$, and in particular, $\dim(I) = \dim(I\setminus(\pi|_C \circ \phi_1)^{-1}(B))$, but $I\setminus(\pi|_C \circ \phi_1)^{-1}(B)$ is not open in $I$, contradicting the claim.

$(1)\Rightarrow(4)$.  Let $T \subseteq S$ be a definable subsemigroup. By Corollary \ref{Forma dos subsemigroups nos completamente simples}, there exists $J \subseteq I$, $\Gamma \subseteq \Lambda$ and a definable subgroup $W \subseteq G$ such that $T = \M J W \Gamma {P|_{\Gamma \times J}}$. By Corollary 2.8 of \cite{pillay_groups_1988}, $W$ is closed, and as $I$ and $\Lambda$ have the discrete topology, $J$ and $\Gamma$ are closed. As a consequence, $T = J \times W \times \Gamma$ is closed.

$(4) \Rightarrow(1)$. Again, we will only show that $I$ is finite as the proof for $\Lambda$ is analogous. Let $x \in I$ and set $I_x = I\setminus \{x\}$. Let $T = \M{I_x}{G}{\Lambda}{P|_{\Lambda \times I_x}}$ be a definable subsemigroup of $S$. By our hypothesis, $T$ is closed, and thus definably compact. Projecting into $I$, we conclude that $I_x$ is definably compact, and thus closed. This means that $\{x\}$ is open in $I$, and as $x \in I$ was arbitrary, we conclude that $I$ is equipped with the discrete topology. Any definable space has only finitely many definably connected components, meaning that $I$ has to be finite. 
\end{cproof}

\textbf{Remark:}
Note that in Proposition \ref{H finite semigroup equivalencia}, definable compactness is necessary. Consider for example the structure $\cM = (\R, <, +, \cdot, 0, 1)$. Then $(\R,+)$ is a definable completely simple semigroup isomorphic to $\M I \R \Lambda P$, where $I = \Lambda = \{0\}$ and $P(0,0) = 0$.

Then $(0,\infty) \supset (1,\infty) \supset \ldots  \supset (n,\infty) \supset\ldots $ is a counterexample to  (2) and $(4)$. Additionally, $[0,\infty)$ is a definable subsemigroup with the same dimension as $\R$ that is not open, thus being a counterexample to (3).

\begin{example}
    Fix the o-minimal structure $\cM = (\R,+,\cdot,0,1)$. Consider the definable semigroup $S = [0,1]^2 \subseteq \R^2$ where the product is given by $$(a,b) * (c,d) = (ac,b),$$
    which is definable. Additionally, with the topology induced by $\R^2$, $S$ is a definably compact semigroup (with a definable manifold structure).

    By Proposition \ref{def compact semigroup has idempotent}, $S$ has at least one idempotent, and indeed, it is easy to verify that $$E(S) = \{0,1\}\times [0,1].$$
    By Proposition \ref{unique minimal ideal}, $S$ has a minimal ideal, which is $$K = \{0\} \times [0,1].$$ It is easy to see that $K$ is an ideal, and for any ideal $I$, for any $(a,b) \in I$ and for any $(0,x) \in K$, we have that $(0,x)(a,b) = (0,x) \in I$, which implies that $K \subseteq I$.

     $K$ has an infinite number of idempotents, meaning $K$ is not $\Hg$-finite. By Corollary \ref{DCC equiv}, $S$ does not have the DCC. Indeed $$K=\{0\}\times [0,1] \supset \{0\}\times [0,1/2] \supset \{0\}\times [0,1/3] \supset \ldots$$
    is an example of an infinite strictly descending chain of definable subsemigroups.

    By Corollary \ref{existem ideais minimais}, $S$ also has minimal left and right ideals. Moreover, for any $e \in E(K) = K$, $Se$ is a minimal left ideal and $eS$ is a minimal right ideal (see \cite{abodayeh_compact_1997} for more details), so we conclude that all minimal right ideals are of the form $$\text{$\{(0,x)\}$ for $x \in [0,1]$}$$
    and the only minimal left ideal is $$\{0\}\times [0,1] = K.$$
\end{example}

We would like to conclude with the following remarks:

Another example where the similarities between definable groups and definable semigroups break is with regards to the following property: given a definable group $G$ and a subset $A \subseteq G$, then there exists a least definable subgroup $\langle A\rangle_{\text{def}}$ that contains $A$. This property comes from the fact that any definable group in an o-minimal structure satisfies the descending chain condition (see Fact \ref{facts about definable groups}). However, if we consider instead a definable semigroup $S$ and a subset $A\subseteq S$, there may not be a least definable subsemigroup $\langle A\rangle _{\text{def}}$ that contains $A$.
To see this, consider $\cM$ to be the ordered field of real numbers, and $S = [0,1]$ with multiplication and the usual topology. For any $\varepsilon\in (0,1)$, consider the set $$S(\varepsilon) = [0,\varepsilon)\cup \{2^{-n} : n \geq 1\}.$$ Then $\{S(\varepsilon) : \varepsilon\in (0,1)\}$ is a strictly decreasing family of subsemigroups containing $2^{-1}$. Furthermore, by o-minimality, any definable subsemigroup containing $1/2$ must contain some $S(\varepsilon)$. As such, $\langle 2^{-1}\rangle_{\text{def}}$ does not exist.

On the other hand, definable compactness seems, so far, to be a good generalization of the notion of compactness for studying definable semigroups, as foundational results in the theory of compact semigroups carry over to definably compact semigroups. However, we would like to pose the following questions that, at this time, we are not sure how to answer.  We say that a semigroup $S$ is \emph{right stable} if, for any $s,x \in S$, $$s \mathrel{\Jg} sx \Rightarrow s \mathrel{\Rg} sx$$
and $S$ is said to be \emph{left stable} if, for any $s,x \in S$, $$s \mathrel{\Jg} xs \Rightarrow s \mathrel{\Lg} xs.$$
Furthermore, we say that $S$ is \emph{stable} if it is right and left stable.

It is known that any compact semigroup is stable, which raises the following natural question:

\begin{question}
    Is it true that any definably compact definable semigroup is also stable?
\end{question}
Another interesting thing to note is that, contrary to the results about compact semigroups, results known to be true about definable groups do not, in general, carry over to this new context, as seen in Proposition \ref{H finite semigroup equivalencia}.

One natural direction for further investigation is whether definable semigroups, like definable groups, admit a definable manifold structure that makes them into topological semigroups.

This leads to the following question:

\begin{question} \label{question definable manifold}
    Given an o-minimal structure $\cM$, and a definable semigroup $S$, does there always exist a definable manifold structure on $S$ that makes it a topological semigroup?
\end{question}

If the answer is positive, in general, said definable manifold structure will not be unique, unlike what happens with Pillay's topology on any definable group. To see this, let $S$ be any definable set in $\cM$, and fix any element of $S$ which we will denote by $0$. The operation $s*t = 0$ for any $s,t \in S$ is associative and definable, meaning that $(S,*)$ is a definable semigroup. Given any subset $X \subseteq S$, then the pre-image of $X$ under $*$ is given by
$$*^{-1}(X) =\begin{cases}
    \emptyset & \text{if $0 \not \in X$}\\
    S \times S & \text{if $0 \in X$}
\end{cases},$$
which implies that any topology on $S$ makes $*$ continuous, and in particular, any definable manifold structure in $S$ will make $S$ a topological semigroup.

\section{Acknowledgments}
This research was conducted in partial fulfillment of the requirements for the Master's degree in Mathematics at the University of Porto, under the supervision of Dr. Mário Edmundo (University of Lisbon) and co-supervision of Dr. Jorge Almeida (University of Porto), whose guidance and support are gratefully acknowledged.

The author is also grateful to the anonymous referee for their careful reading of the manuscript and for many helpful comments and suggestions, which led to significant improvements in the exposition and in the overall presentation of the paper.

\bibliographystyle{plain}
\bibliography{references}

\end{document}